\def\autori{P. DALL'AGLIO and C. LEONE}
\def\titolo{Obstacles problems with measure data} 




\font\sixrm=cmr6
\newcount\tagno \tagno=0                        
\newcount\thmno \thmno=0                        
\newcount\bibno \bibno=0                        
\newcount\chapno\chapno=0                       
\newcount\verno            
\newif\ifproofmode
\proofmodetrue
\newif\ifwanted
\wantedfalse
\newif\ifindexed
\indexedfalse

\def\ifundefined#1{\expandafter\ifx\csname+#1\endcsname\relax}

\def\Wanted#1{\ifundefined{#1} \wantedtrue 
\immediate\write0{Wanted 
#1 
\the\chapno.\the\thmno}\fi}

\def\Increase#1{{\global\advance#1 by 1}}

\def\Assign#1#2{\immediate
\write1{\noexpand\expandafter\noexpand\def
 \noexpand\csname+#1\endcsname{#2}}\relax
 \global\expandafter\edef\csname+#1\endcsname{#2}}

\def\pAssign#1#2{\write1{\noexpand\expandafter\noexpand\def
 \noexpand\csname+#1\endcsname{#2}}}

\def\lPut#1{\ifproofmode\llap{\hbox{\sixrm #1\ \ \ }}\fi}
\def\rPut#1{\ifproofmode$^{\hbox{\sixrm #1}}$\fi}



\def\chp#1{\global\tagno=0\global\thmno=0\Increase\chapno
\Assign{#1}
{\the\chapno}{\lPut{#1}\the\chapno}}


\def\thm#1{\Increase\thmno 
\Assign{#1}
{\the\chapno.\the\thmno}\the\chapno.\the\thmno\rPut{#1}}


\def\frm#1{\Increase\tagno
\Assign{#1}{\the\chapno.\the\tagno}\lPut{#1}
{\the\chapno.\the\tagno}}


\def\bib#1{\Increase\bibno
\Assign{#1}{\the\bibno}\lPut{#1}{\the\bibno}}


\def\pgp#1{\pAssign{#1/}{\the\pageno}}


\def\ix#1#2#3{\pAssign{#2}{\the\pageno}
\immediate\write#1{\noexpand\idxitem{#3}
{\noexpand\csname+#2\endcsname}}}

 
\def\rf#1{\Wanted{#1}\csname+#1\endcsname\relax\rPut {#1}}


\def\rfp#1{\Wanted{#1}\csname+#1/\endcsname\relax\rPut{#1}}

\input \jobname.auxi
\Increase\verno
\immediate\openout1=\jobname.auxi

\immediate\write1{\noexpand\verno=\the\verno}

\ifindexed
\immediate\openout2=\jobname.idx
\immediate\openout3=\jobname.sym 
\fi


\font\twelverm=cmr12
\font\twelvei=cmmi12
\font\twelvesy=cmsy10
\font\twelvebf=cmbx12
\font\twelvett=cmtt12
\font\twelveit=cmti12
\font\twelvesl=cmsl12

\font\ninerm=cmr9
\font\ninei=cmmi9
\font\ninesy=cmsy9
\font\ninebf=cmbx9
\font\ninett=cmtt9
\font\nineit=cmti9
\font\ninesl=cmsl9

\font\eightrm=cmr8
\font\eighti=cmmi8
\font\eightsy=cmsy8
\font\eightbf=cmbx8
\font\eighttt=cmtt8
\font\eightit=cmti8
\font\eightsl=cmsl8

\font\sixrm=cmr6
\font\sixi=cmmi6
\font\sixsy=cmsy6
\font\sixbf=cmbx6

\catcode`@=11 
\newskip\ttglue

\def\twelvepoint{\def\rm{\fam0\twelverm}
\textfont0=\twelverm  \scriptfont0=\ninerm  
\scriptscriptfont0=\sevenrm
\textfont1=\twelvei  \scriptfont1=\ninei  \scriptscriptfont1=\seveni
\textfont2=\twelvesy  \scriptfont2=\ninesy  
\scriptscriptfont2=\sevensy
\textfont3=\tenex  \scriptfont3=\tenex  \scriptscriptfont3=\tenex
\textfont\itfam=\twelveit  \def\it{\fam\itfam\twelveit}%
\textfont\slfam=\twelvesl  \def\sl{\fam\slfam\twelvesl}%
\textfont\ttfam=\twelvett  \def\tt{\fam\ttfam\twelvett}%
\textfont\bffam=\twelvebf  \scriptfont\bffam=\ninebf
\scriptscriptfont\bffam=\sevenbf  \def\bf{\fam\bffam\twelvebf}%
\tt  \ttglue=.5em plus.25em minus.15em
\normalbaselineskip=15pt
\setbox\strutbox=\hbox{\vrule height10pt depth5pt width0pt}%
\let\sc=\tenrm  \let\big=\twelvebig  \normalbaselines\rm}

\def\tenpoint{\def\rm{\fam0\tenrm}
\textfont0=\tenrm  \scriptfont0=\sevenrm  \scriptscriptfont0=\fiverm
\textfont1=\teni  \scriptfont1=\seveni  \scriptscriptfont1=\fivei
\textfont2=\tensy  \scriptfont2=\sevensy  \scriptscriptfont2=\fivesy
\textfont3=\tenex  \scriptfont3=\tenex  \scriptscriptfont3=\tenex
\textfont\itfam=\tenit  \def\it{\fam\itfam\tenit}%
\textfont\slfam=\tensl  \def\sl{\fam\slfam\tensl}%
\textfont\ttfam=\tentt  \def\tt{\fam\ttfam\tentt}%
\textfont\bffam=\tenbf  \scriptfont\bffam=\sevenbf
\scriptscriptfont\bffam=\fivebf  \def\bf{\fam\bffam\tenbf}%
\tt  \ttglue=.5em plus.25em minus.15em
\normalbaselineskip=12pt
\setbox\strutbox=\hbox{\vrule height8.5pt depth3.5pt width0pt}%
\let\sc=\eightrm  \let\big=\tenbig  \normalbaselines\rm}

\def\ninepoint{\def\rm{\fam0\ninerm}
\textfont0=\ninerm  \scriptfont0=\sixrm  \scriptscriptfont0=\fiverm
\textfont1=\ninei  \scriptfont1=\sixi  \scriptscriptfont1=\fivei
\textfont2=\ninesy  \scriptfont2=\sixsy  \scriptscriptfont2=\fivesy
\textfont3=\tenex  \scriptfont3=\tenex  \scriptscriptfont3=\tenex
\textfont\itfam=\nineit  \def\it{\fam\itfam\nineit}%
\textfont\slfam=\ninesl  \def\sl{\fam\slfam\ninesl}%
\textfont\ttfam=\ninett  \def\tt{\fam\ttfam\ninett}%
\textfont\bffam=\ninebf  \scriptfont\bffam=\sixbf
\scriptscriptfont\bffam=\fivebf  \def\bf{\fam\bffam\ninebf}%
\tt  \ttglue=.5em plus.25em minus.15em
\normalbaselineskip=11pt
\setbox\strutbox=\hbox{\vrule height8pt depth3pt width0pt}%
\let\sc=\sevenrm  \let\big=\ninebig  \normalbaselines\rm}

\def\eightpoint{\def\rm{\fam0\eightrm}
\textfont0=\eightrm  \scriptfont0=\sixrm  \scriptscriptfont0=\fiverm
\textfont1=\eighti  \scriptfont1=\sixi  \scriptscriptfont1=\fivei
\textfont2=\eightsy  \scriptfont2=\sixsy  \scriptscriptfont2=\fivesy
\textfont3=\tenex  \scriptfont3=\tenex  \scriptscriptfont3=\tenex
\textfont\itfam=\eightit  \def\it{\fam\itfam\eightit}%
\textfont\slfam=\eightsl  \def\sl{\fam\slfam\eightsl}%
\textfont\ttfam=\eighttt  \def\tt{\fam\ttfam\eighttt}%
\textfont\bffam=\eightbf  \scriptfont\bffam=\sixbf
\scriptscriptfont\bffam=\fivebf  \def\bf{\fam\bffam\eightbf}%
\tt  \ttglue=.5em plus.25em minus.15em
\normalbaselineskip=9pt
\setbox\strutbox=\hbox{\vrule height7pt depth2pt width0pt}%
\let\sc=\sixrm  \let\big=\eightbig  \normalbaselines\rm}

\def\twelvebig#1{{\hbox{$\textfont0=\twelverm\textfont2=\twelvesy
	\left#1\vbox to10pt{}\right.\n@space$}}}
\def\tenbig#1{{\hbox{$\left#1\vbox to8.5pt{}\right.\n@space$}}}
\def\ninebig#1{{\hbox{$\textfont0=\tenrm\textfont2=\tensy
	\left#1\vbox to7.25pt{}\right.\n@space$}}}
\def\eightbig#1{{\hbox{$\textfont0=\ninerm\textfont2=\ninesy
	\left#1\vbox to6.5pt{}\right.\n@space$}}}
 
\def\displayliness#1{\null\,\vcenter{\openup1\jot \m@th
  \ialign{\strut\hfil$\displaystyle{##}$\hfil
    \crcr#1\crcr}}\,}
	       
\def\displaylinesno#1{\displ@y \tabskip=\centering
   \halign to\displaywidth{ \hfil$\@lign \displaystyle{##}$ \hfil
	\tabskip=\centering
     &\llap{$\@lign##$}\tabskip=0pt \crcr#1\crcr}}
		     
\def\ldisplaylinesno#1{\displ@y \tabskip=\centering
   \halign to\displaywidth{ \hfil$\@lign \displaystyle{##}$\hfil
	\tabskip=\centering
     &\kern-\displaywidth
     \rlap{$\@lign##$}\tabskip=\displaywidth \crcr#1\crcr}}

\catcode`@=12 

\def\parag#1#2{\goodbreak\bigskip\bigskip\noindent
		   {\bf #1.\ \ #2}
		   \nobreak\bigskip} 
\def\intro#1{\goodbreak\bigskip\bigskip\goodbreak\noindent
		   {\bf #1}\nobreak\bigskip\nobreak}
\long\def\th#1#2{\goodbreak\bigskip\noindent
		{\bf Theorem #1.\ \ \it #2}}
\long\def\lemma#1#2{\goodbreak\bigskip\noindent
		{\bf Lemma #1.\ \ \it #2}}
\long\def\prop#1#2{\goodbreak\bigskip\noindent
		  {\bf Proposition #1.\ \ \it #2}}

\long\def\defin#1#2{\goodbreak\bigskip\noindent
		  {\bf Definition #1.\ \ \rm #2}}
\long\def\rem#1#2{\goodbreak\bigskip\noindent
		 {\bf Remark #1.\ \ \rm #2}}
\long\def\ex#1#2{\goodbreak\bigskip\noindent
		 {\bf Example #1.\ \ \rm #2}}

\def\negbigskip{\vskip-.4cm}

\def\proof{\vskip.4cm\noindent{\it Proof.\ \ }}

\def\sqr#1#2{\vbox{
   \hrule height .#2pt 
   \hbox{\vrule width .#2pt height #1pt \kern #1pt 
      \vrule width .#2pt}
   \hrule height .#2pt }}
\def\square{\sqr74}

\def\endproof{{\unskip\nobreak\hfill \penalty50
\hskip2em\hbox{}\nobreak\hfill $\square$ \goodbreak
\parfillskip=0pt  \finalhyphendemerits=0}}

\mathchardef\emptyset="001F
\mathchardef\hyphen="002D


\def\wto{\rightharpoonup}


\def\rightheadline{\eightpoint\hfil\titolo
\hfil\tenrm\folio} 
\def\leftheadline{\tenrm\folio\hfil\eightpoint
\autori \hfil}
\def\zeroheadline{\hfill} 
%
\headline={\ifnum\pageno=0 \zeroheadline
\else\ifodd\pageno\rightheadline
\else\leftheadline\fi\fi}

\nopagenumbers
\magnification=1200
\baselineskip=15pt
\hfuzz=2pt
\parindent=2em
\mathsurround=1pt
\tolerance=1000

\pageno=0
\hsize 14truecm
\vsize 25truecm
\hoffset=0.8truecm
\voffset=-1.55truecm

\null
\vskip 2.8truecm
{\twelvepoint
\baselineskip=1.7\baselineskip
\centerline{\bf OBSTACLE PROBLEMS WITH}
\centerline{\bf MEASURE DATA}
}
\vskip2truecm
\centerline{Paolo DALL'AGLIO}
\smallskip
\centerline{Chiara LEONE}
\vfil

{\eightpoint
\baselineskip=1.2\baselineskip
\centerline{\bf Abstract}
\medskip
\noindent 
We give a definition for Obstacle Problems with measure data and general 
obstacles. For such problems we prove existence and uniqueness of 
solutions and consistency with the classical theory of Variational 
Inequalities. Continuous dependence with respect to data is discussed.
\bigskip
}

\vfil
\vskip 1truecm 
\centerline {Ref. S.I.S.S.A. 
147/97/M (November 97)}
\vskip 1truecm 
\eject

%
%
\topskip=25pt 
\vsize 22.5truecm
\hsize 16.2truecm
\hoffset=0truecm
\voffset=0.5truecm

\proofmodefalse 

\def\R{\rm I\! R}  
\def\Luno{{\rm L}^1(\Omega)}
\def\Linf{{\rm L}^\infty(\Omega)}
\def\Hunozero{{\rm H}^1_0(\Omega)}
\def\Hduale{{\rm H}^{\hbox{\kern 1truept \rm -}\kern -1truept1}\kern
	-1truept(\Omega)}
\def\Wunoqzero{W^{1,q}_0(\Omega)}
\def\Mb{{\cal M}_b(\Omega)}
\def\Mbp{{\cal M}^+_b(\Omega)}
\def\Mbo{{\cal M}^0_b(\Omega)}
\def\Inters{\Mb\cap\Hduale}
\def\K{{\rm K}_\psi(\Omega)}
\def\intl{\int\limits}
\def\inters{\K\cap W^{1,q}_0(\Omega)}
\def\dualita#1#2{\langle#1,#2\rangle}
\def\Fpsimu{{\cal F}_\psi(\mu)}

\def\wto{\rightharpoonup}

\def\um{u_{\mu}}
\def\un{u_{\nu}}
\def\ul{u_{\lambda}}
\def\weak{\hbox{ \ weakly in }\Wunoqzero}
\def\weak2{\hbox{ \ weakly in }\Hunozero}
\def\weaks{\hbox{ \ $\ast$-weakly in }\Mb}
\def\strongq{\hbox{ strongly in }\Wunoqzero}
\def\ae{\hbox{ \ a.e. in }\Omega}
\def\qe{\hbox{ \ q.e. in }\Omega}
\def\vuoto{{\scriptstyle\bigcirc\!\!\!}\!\lower.15ex\hbox{/}}
\def\f{\varphi}
\def\e{\varepsilon}
\def\fe{\varphi_\varepsilon}
\def\crr{\cr\noalign{\medskip}}

\parag{\chp{int}}{Introduction}

In this paper we consider the obstacle problem with measure data for a 
linear differential operator ${\cal A}$, for which we prove existence and 
uniqueness of solutions together with some stability results.

Consider first the objects that won't change throughout the work. 

Let $\Omega$ be a regular bounded open subset of $\R^N$ (for the notion
of regularity see Definition~\rf{regolarita}).

Let ${\cal A}(u)=-{\rm div}(A(x)\nabla u)$ be a linear elliptic operator
with coefficients in $L^\infty(\Omega)$, that is $A(x)=((a_{ij}(x)))$ is 
an $N\times N$ matrix such that
$$
a_{ij}\in L^\infty(\Omega)\hbox{\ \ and \ }\sum a_{ij}(x)\xi_i\xi_j\geq\gamma
|\xi|^2,\quad\forall\xi\in\R^N, \ae.
$$

\bigskip
We want to consider the obstacle problem also in the case of thin
obstacles, so we will need the techniques of capacity theory. For this
theory we refer, for instance, to [\rf{HEI-KIL}].

We recall very briefly that, given a set $E\subseteq\Omega$, its
capacity with respect to $\Omega$ is given by  
$$
{\rm cap}(E,\Omega)=\inf\{\|z\|^2_{\Hunozero}\,:\,z\in\Hunozero, z\geq
1\hbox{ a.e. in a neighbourhood of }E\}.
$$
When the ambient set $\Omega$ is clear from the context we will write
cap$(E)$.
        
A property holds quasi everywhere (abbreviated as q.e.) when it holds up
to sets of capacity zero.

A function $v:\Omega\to{\overline{\R}}$ is quasi continuous (resp.
quasi upper semicontinuous) if, for every $\e>0$ there exists a set
$E$ such that ${\rm cap}(E)<\e$ and $v_{|_{\Omega\setminus E}}$ is
continuous
(resp. upper semicontinuous) in
$\Omega\setminus E$.

We recall also that if $u$ and $v$ are quasi continuous functions and
$u\leq v$ a.e. in $\Omega$ then also $u\leq v$ q.e. in $\Omega$.

A function $u\in\Hunozero$ always has a quasi continuous representative,
that is there
exists a quasi continuous function $\tilde u$ which equals $u$ a.e. in
$\Omega$. We shall always identify $u$ with its quasi continuous
representative. We also have that if $u_n\to u$ strongly in $\Hunozero$
then there
exists a subsequence which converges quasi everywhere.


\bigskip

Consider the function $\psi:\Omega\to{\overline{\R}}$, and let the convex 
set be
$$\K:=\{z \hbox{ quasi continuous }\ :\ z\geq\psi\qe\}.$$

Without loss of generality we may always assume that $\psi$ is quasi
upper semicontinuous thanks to Proposition
1.5 in [\rf{DAL-4}].

\bigskip

In their natural setting, obstacle problems are part of 
the theory of Variational
Inequalities (for which we refer to well known books such as
[\rf{KIN-STA}] and [\rf{TRO}]).

For any datum $F\in\Hduale$ the
Variational Inequality with obstacle $\psi$
$$
\cases{\dualita{{\cal A} u}{v-u}\geq\dualita{F}{v-u}
		\quad\forall v\in\K\cap\Hunozero\crr
	u\in\K\cap\Hunozero\cr}\eqno(\frm{vi})
$$
(which, for simplicity, will be indicated by
$
VI(F,\psi)
$),
has a unique solution, whenever the set $\K\cap \Hunozero$ is nonempty, 
that is ensured by the condition
$$
\exists z\in\Hunozero\ :\ z\geq\psi\qe.\eqno(\frm{zeta})
$$

In this frame, among all classical results, we recall 
that the 
solution of $VI(F,\psi)$ is also characterized as 
the smallest function $u\in\Hunozero$ such that
$$
\cases{
{\cal A}u-F\geq 0
	\hbox{ in }{\cal D}'(\Omega)\crr
u\geq\psi\qe.\cr}\eqno(\frm{caratt})
$$
Then $\lambda:={\cal A}u-F$ is a positive measure, and will be called 
the obstacle reaction associated with $u$.
\bigskip

In order to study the case of right-hand side measure, 
we recall that, already in the case of equations, the term $\langle 
\mu,u\rangle$ has not always meaning when $\mu$ is a measure and 
$u\in\Wunoqzero,\,q<N$. 
Hence the classical formulation of the variational inequality fails.

Also the use of the characterization (\rf{caratt}) 
to define the Obstacle 
Problem with measure data is not possible because another problem 
arises: a famous example 
by J. Serrin (see [\rf{SER}] and, for more details, [\rf{PRI}]) shows that 
the homogeneous equation
$$
\cases{
	{\cal A}u=0\hbox{ in }{\cal D}'(\Omega)\crr
	u=0\hbox { on }\partial\Omega\cr}
$$
has a nontrivial solution $v$ which does not belong to $\Hunozero$. Here 
${\cal A}$ is a particular linear elliptic operator with discontinuous 
coefficients. The function $u=0$ is obviously the unique solution in 
$\Hunozero$.

So (\rf{caratt}) in general does not determine the solution of the 
Obstacle Problem: indeed, with such ${\cal A}$, if we choose $\psi\equiv 
-\infty$, and if $u$ were the minimal supersolution, then we would have 
$u\leq u+ tv\ae$ for any $t$ in $\R$, which is a contradiction.

G. Stampacchia overcame this difficulty for equations, using a wider 
class of test functions, and gave in [\rf{STA}] the 
following definition, which uses regularity and duality arguments.

For this theory we need to assume that the boundary $\partial\Omega$ has
the following property, which is satisfied in particular when
$\partial\Omega$ is Lipschitz.

\defin{\thm{regolarita}}{}We say that a bounded open subset of $\R^N$ is
regular if there exists a constant $\alpha\in (0,1)$ such that, for every
$x_0\in\partial\Omega$ and for all $\rho$ small, we have
$$
|B_\rho(x_0)\setminus\Omega|>\alpha\,|B_\rho(x_0)|.
$$

The theory by Stampacchia actually works under slightly
weaker but more complicated assumptions as said in [\rf{STA}] (see
Definition~6.2).

Let  $\mu\in\Mb$, where $\Mb$ is the space of bounded Radon measures, viewed 
as the dual of the Banach space $C_0(\Omega)$ of continuous functions that 
are zero on the boundary. 

\defin{\thm{stamp}}{}
A function $\um\in\Luno$ is a solution
in the sense of Stampacchia (also called 
solution by duality) of the
equation
$$
\cases{{\cal A}\um=\mu\ \ \hbox{ in }\Omega\crr
	\um=0\hbox{ on }\partial\Omega,\cr}
$$
if
$$
\intl_\Omega \um g\,dx\,=\,\intl_\Omega u^*_g\,d\mu,\quad\forall g\in\Linf,
$$
where $u^*_g$ is the solution of
$$
\cases{{\cal A}^*u^*_g=g\quad \hbox{ in }\Hduale\crr
	u^*_g\in\Hunozero\cr}
$$
and ${\cal A}^*$ is the adjoint of ${\cal A}$. 
\bigskip

Stampacchia proved that a solution $\um$ exists and is unique and
belongs to $\Wunoqzero$, where $q$ is any exponent satisfying
$1<q<\displaystyle{ N\over{N-1}}$.  He proved also that,
if the datum $\mu$ is more regular, namely in $\Inters$, then the solution
coincides with the variational one.
It is also known that $T_k(\um)$ belongs to $\Hunozero$, where $T_k$
denotes the usual truncation function defined by
$$
T_k(s):=(-k)\vee(s\wedge k).
$$
From this it follows that every Stampacchia solution $\um$ has a quasi
continuous representative; in the rest of the paper we shall always
identify $\um$ with its quasi continuous representative.

Moreover when the data converge
$*$-weakly in $\Mb$, also the solutions converge strongly in $\Wunoqzero$ 
and their truncates weakly in $\Hunozero$.

We will use the following notation: $\um$ denotes 
the solution of the equation 
$$
\cases{{\cal A}\um=\mu\quad\hbox{ in }\Omega\crr
	\um=0\qquad\hbox{ on }\partial\Omega,\cr}
$$
when $\mu$ is either a measure in $\Mb$ or an element of $\Hduale$. In the 
first case we refer to the definition by G. Stampacchia, in the latter 
to the usual variational one.

\bigskip

Following these ideas we give a formulation 
for Obstacle Problems which involves this type of solutions.

\defin{\thm{disvar}}{}We say that the function 
$u\in\inters,\,1<q<{\displaystyle{N\over {N-1}}}$ is a solution of the
Obstacle Problem with datum $\mu$ and obstacle $\psi$
if 
\item{1.} there exists a positive bounded measure $\lambda\in\Mbp$ such
that
$$
u=\um+\ul;
$$
\item{2.} for any $\nu\in\Mbp$, such that $v=\um+\un$ belongs to $\K$, 
we have
$$
u\leq v\qe.
$$

Also here the positive measure $\lambda$, which is uniquely defined,
will be called the  
obstacle reaction relative to $u$.
This problem will be shortly indicated by $OP(\mu,\psi)$.

\bigskip

To show that for any datum $\mu$ there exists one and only one solution,
we introduce the set 
$$
\Fpsimu:=\left\{v\in\inters\ :\ \exists\nu\in\Mbp\ 
	\hbox{s.t.}\ v=\um+\un\right\}.
$$
We will prove that $\Fpsimu$ has a minimum element, that is a
function $u\in\Fpsimu$ such that $u\leq v\qe$ for any other function 
$v\in\Fpsimu$.
This is clearly the solution of the Obstacle Problem according to
the Definition \rf{disvar}.
If this solution exists it is obviously unique.

Hypothesis (\rf{zeta}) does not ensure that 
$\Fpsimu$ be nonempty. The 
minimal hypothesis, instead of (\rf{zeta}), 
will be 
$$
\exists \rho\in\Mb\ :\ u_\rho\geq\psi\hbox{ q.e. in
}\Omega;\eqno(\frm{ipomin})
$$
so the set $\Fpsimu$ is nonempty for 
every $\mu\in\Mb$, because it contains the function $u_{\mu^+}+u_\rho$.

The proof of existence will be first worked out for the case of 
a negative obstacle (Section \rf{psinegativo}): this
is based on an approximation
technique. The obstacle reactions associated with the 
solutions for regular data are shown to satisfy 
an estimate on the masses, which allows to pass to the limit and 
obtain the solution in the general case. Then the proof is easily extended 
to the case of general obstacle (Section \rf{finedimondo})

In Section \rf{nuovavita} we will give some stability results, and in Section 
\rf{ipotesi} we will show that the classical 
solution to the Obstacle Problem (equation (\rf{vi}))  
coincides with the new one (Definition \rf{disvar})
when both make sense.

Moreover we will show that this solution coincides with the one given
in the wider setting of nonlinear monotone operators, but in the 
case of datum in $\Luno$, by L. 
Boccardo and T. Gallou\"et in [\rf{BOC-GAL}], and by L. Boccardo an G.R. 
Cirmi in [\rf{BOC-CIR}] and [\rf{BOC-CIR2}].

Section \rf{approx} provides a characterization of the solution in terms 
of appoximating sequences of solutions of Variational Inequalities, and 
in Section \rf{emmezero} we study some properties of the solution of the 
Obstacle Problem for the class of Radon measures, that do not 
charge the sets of zero 2-Capacity.

\parag{\chp{psinegativo}}{Nonpositive obstacles}

Throughout this chapter we assume the obstacle to be nonpositive. In this 
frame both hypotheses (\rf{ipomin}) and (\rf{zeta}) are trivially 
satisfied.

We begin with a preparatory result which will be proved in two steps.

\lemma{\thm{mupiuemumeno}}{Let $\psi\leq0$ an let $\mu\in\Inters$ such 
that $\mu^+$ 
and $\mu^-$ belong to $\Hduale$. Let $u$ be the solution of  
$VI(\mu,\psi)$ and $\lambda$ the obstacle reaction associated with  $u$. 
Then 
$$
||\lambda||_{\Mb}\,\leq\,||\mu^-||_{\Mb}.
$$}
\negbigskip
\negbigskip
\proof{}Observe that the function $u_{\mu^+}$ is positive and hence 
greater than or equal to $\psi$, belongs to $\Hunozero$, and
$$
{\cal A}u_{\mu^+}-\mu\geq0\quad\hbox{ in }{\cal D}'(\Omega).
$$
By (\rf{caratt}) we
have
$$
u=\um+\ul\leq u_{\mu^+}\qe,
$$
and, by linearity,
$$
\ul\leq u_{\mu^-}\qe.\eqno(\frm{mumeno})
$$
We will prove that this implies 
$$
\lambda(\Omega)\leq\mu^-(\Omega)\eqno(\frm{tesi})
$$
which is equivalent to the thesis.

To prove (\rf{tesi}) we note that, thanks to (\rf{mumeno})
$$
\intl_\Omega w\,d\mu^-=\dualita{{\cal A}^*w}{u_{\mu^-}}\geq
\dualita{{\cal A}^*w}{\ul}=\intl_\Omega w\, 
d\lambda,\eqno(\frm{quattro})
$$
for every $w\in\Hunozero$, such that ${\cal A}^*w\geq 0$ in ${\cal 
D}'(\Omega)$.

It is now easy to find a sequence $\{w_n\}$ in $\Hunozero$ such that
$w_n\nearrow 1$ and ${\cal A}^*w_n\geq 0$ in ${\cal
D}'(\Omega)$. For instance, one can choose as $w_n$ the ${\cal 
A}^*$-capacitary potential (see [\rf{HEI-KIL}], chapter 9) of $J_n$, where 
$J_n$ is an invading family of compact subsets of $\Omega$.

Passing to the limit in (\rf{quattro}), as $n\to\infty$, we obtain 
(\rf{tesi}).
\endproof

\th{\thm{masse}}{Let $\psi\leq0$ and $\mu\in\Inters$. Let $u$ be the 
solution of $VI(\mu,\psi)$ and let $\lambda$ be the obstacle reaction 
relative to $u$. Then $$
||\lambda||_{\Mb}\leq||\mu^-||_{\Mb}.\eqno(\frm{Lewy})
$$}
\negbigskip

\negbigskip
\proof Thanks to Lemma 3.3 in [\rf{DM-MAL}] there exists a sequence of
smooth functions $f_n$ such that
$$
||f_n-\mu||_{\Hduale}\leq{1\over n}\quad\hbox{  and }\quad
||f_n||_{\Luno}\leq||\mu||_{\Mb}.
$$
Thanks to the next Lemma, the sequence $f_n$ satisfies
$$
f^{\pm}_n\wto\mu^{\pm}\weaks\hbox{ and }
||f^{\pm}_n||_{\Luno}\to||\mu^{\pm}||_{\Mb}.
$$

Let $u_n$ and $u$ be the solutions of
$VI(f_n,\psi)$ and $VI(\mu,\psi)$, respectively. 
We know from the general 
theory (see, for instance, [\rf{KIN-STA}]) that $u_n\to u$ in 
$\Hunozero$. So the measures $\lambda_n$ 
and $\lambda$ associated with $u_n$ and $u$, respectively, satisfy
$$
\lambda_n\to\lambda\hbox{ in }\Hduale,
$$
$$
||\lambda_n||_{\Mb}\leq||f^-_n||_{\Luno}.
$$
So $\lambda_n\wto\lambda$ in$\weaks$, and we get the inequality
(\rf{Lewy}).
\endproof

The following lemma is quite simple, but is proved here for the sake 
of completeness.

\lemma{\thm{piumeno}}{Let $\mu_n$ and $\mu$ be measures in $\Mb$ such 
that 
$$
\mu_n\wto\mu\weaks\hbox{ and }||\mu_n||_{\Mb}\to||\mu||_{\Mb}
$$
then
$$
\mu_n^+\wto\mu^+\hbox{ and }\mu^-_n\wto\mu^-\weaks,
$$
and
$$
\|\mu_n^+\|_{\Mb}\to\|\mu^+\|_{\Mb}\hbox{ and }
\|\mu_n^-\|_{\Mb}\to\|\mu^-\|_{\Mb}.\eqno(\frm{norme})
$$}

\negbigskip

\negbigskip

\proof
Observe that 
$$
||\mu^{\pm}_n||_{\Mb}\leq||\mu_n||_{\Mb}.
$$
so, up to a subsequence, 
$$
\mu^{+}_n\wto\alpha\hbox{ and } \mu^-_n\wto\beta\weaks;
$$
where $\alpha-\beta=\mu$.
Hence, we can compute
$$ 
\eqalign{||\alpha||_{\Mb}+||\beta||_{\Mb}&\leq
		\liminf||\mu^+_n||_{\Mb}+\liminf||\mu^-_n||_{\Mb}\crr
	&\leq\liminf||\mu_n||_{\Mb}=||\mu||_{\Mb};\cr}
$$
from which we easily deduce that $\alpha=\mu^+$,  $\beta=\mu^-$.
Therefore the whole sequences $\mu_n^+$ and $\mu_n^-$ converge to $\mu^+$
and $\mu^-$ respectively. Moreover, as
$$
\eqalign{
\limsup_{n\to+\infty}\|\mu_n^+\|_{\Mb}&+
\liminf_{n\to+\infty}\|\mu_n^-\|_{\Mb}\crr
&\leq
\lim_{n\to+\infty}\|\mu_n\|_{\Mb}=||\mu||_{\Mb}=\|\mu^+\|_{\Mb}+
\|\mu^-\|_{\Mb}\cr}
$$
we obtain easily the first relation in (\rf{norme}). The second one is   
obtained in a similar way.

\endproof

\bigskip

In order to proceed we need to prove that when both the classical 
formulation for the obstacle problem and the new one, given in Definition
\rf{disvar}, make sense
then the solutions, when they exist, are the same. At present we prove it for
a nonpositive obstacle, and we will prove it in the general case in
Section \rf{ipotesi}.

\lemma{\thm{equiv1}}{Let $\mu$ be an element of $\Inters$ and $\psi$ a
nonpositive function; then the solution of $VI(\mu,\psi)$ coincides with 
the solution of $OP(\mu,\psi)$.}
\proof
Let $u$ be the solution of $VI(\mu,\psi)$ and $\lambda$ be the
corresponding obstacle reaction.
Thanks to Theorem~\rf{masse} it is an element of $\Mb$; so
$u\in\Fpsimu$.
Take $v$ an element in $\Fpsimu$, then $v=u_{\mu}+u_{\nu}$, with
$\nu\in\Mbp$, and $v\geq\psi\qe$.

Consider the approximation of
$\nu$, given by ${\cal A}T_k(\un)=:\nu_k$. This is such that $\nu_k\wto
\nu\weaks$ and $\nu_k\in\Mbp\cap\Hduale$ (see [\rf{BOC-MUR}]). Set
$v_k=u_{\mu}+u_{\nu_k}=u_{\mu}+T_k(\un)$. Since trivially
$T_k(\un)\nearrow \un\qe$, we have
$$
v_k\nearrow v\qe.
$$
Denote now the solutions of $VI(\mu,\psi_k)$ by $u_k$, where $\psi_k$ are
the functions defined by
$$
\psi_k:=\psi\wedge v_k.
$$
From $\psi_k\leq\psi_{k+1}\qe$ it easily follows that $u_k\leq
u_{k+1}\qe$. Then there exists a function $u^*$
such that $u_k\nearrow u^*\qe$. 

So, passing to the limit in $u_k\geq\psi_k\qe$ we obtain $u^*\geq
\psi\qe$.

Moreover it is easy to see that $\|u_k\|_{\Hunozero}\leq C$. So, thanks to
Lemma~1.2 in [\rf{DAL-GAR}] we get that $u^*$ is a quasi continuous
function of $\Hunozero$ such that
$$
u_k\wto u^*\weak2.
$$

Moreover it can be easily proved that  $u^*$ is the solution of
$VI(\mu,\psi)$ and by uniqueness $u^*=u\qe$.

Naturally, from the minimality of $u_k$, we deduce
$$
u_k\leq v_k\qe.
$$
so, passing to the limit as $k\to+\infty$ we conclude that $u\leq v\qe$.
Since this is true for every $v\in\Fpsimu$, the function $u$ is the
minimum  in $\Fpsimu$,
i.e. the solution of $OP(\mu,\psi)$.
\endproof 

\bigskip

We are now in position to prove 
that, for every $\mu\in\Mb$ and for every $\psi\leq0$, there exists a 
solution to the Obstacle Problem according to Definition \rf{disvar}.

\th{\thm{esistnegativo}}{Let $\psi\leq0$ and $\mu\in\Mb$.
Then there exists
a solution of $OP(\mu,\psi)$.}

\proof
Consider the function $\um$ and define
$$
{\cal A}(T_k(\um))=:\mu_k.
$$
We know from [\rf{BOC-MUR}] that
$$
\mu_k\wto \mu\weaks
$$
and $\mu_k\in\Hduale$.

Let $u_k$ be the solution of $VI(\mu_k,\psi)$ and denote 
$$
{\cal A}u_k-\mu_k=:\lambda_k,
$$
which we know from Theorem~\rf{masse} to be a measure in $\Mbp$
such that
$$
\|\lambda_k\|_{\Mb}\leq\|\mu_k^-\|_{\Mb}.\eqno(\frm{Lewyk})
$$

Up to a subsequence $\lambda_k\wto\lambda\weaks$, $u_k\to u$ strongly in 
$\Wunoqzero$,
with $u=\um+\ul$, and also $T_h(u_k)\wto T_h(u)$ weakly in
$\Hunozero$, for all
$h>0$.

Now the set
$$
E:=\{v\in\Hunozero\ :\ v\geq T_h(\psi)\qe\}
$$
is closed and convex in $\Hunozero$,
so it is
also weakly closed.
Since, clearly, $T_h(u_k)\geq T_h(\psi)\qe$, passing to the limit as
$k\to+\infty$ we get that also $T_h(u)\in E$, hence $T_h(u)\geq
T_h(\psi)\qe$ for all $h>0$. Passing to the
limit as $h\to+\infty$ we get $u\geq \psi\qe$
In conclusion we deduce $u\in\Fpsimu$.

To show that $u$ is minimal, take $v\in\Fpsimu$ so that $v\geq\psi$ and
$v=\um+\un$.

Let $v_k=u_{\mu_k}+\un$ so that
$v_k=T_k(\um)+\un$ and $v_k\to v$ strongly in $\Wunoqzero$.

Since $\psi\leq0$, we have that $v_k\geq\psi\qe$. 
As $u_k$ is the minimum of ${\cal F}_\psi(\mu_k)$, by Lemma \rf{equiv1}, 
we obtain $u_k\leq v_k\ae$ and in the limit $u\leq v\ae$. Hence $u$
solves $OP(\mu,\psi)$.
\endproof
\bigskip
From formula (\rf{Lewyk}) we see that to extend (\rf{Lewy}) to the case of 
$\mu\in\Mb$ we just need to show that 
$$
||\mu^-_k||_{\Mb}\to||\mu^-||_{\Mb};
$$
this is proved in the following Proposition.

\prop{\thm{Lewy2}}{Let $\psi\leq 0$ and $\mu\in\Mb$. Let $u$ be the 
solution of $OP(\mu,\psi)$ and $\lambda$ the 
corresponding obstacle reaction. Then
$$
||\lambda||_{\Mb}\leq||\mu^-||_{\Mb}.
$$}
\negbigskip
\negbigskip
\proof
What we need is implicit in [\rf{BOC-MUR}]; we recall the main steps of 
that proof, having a closer look to the constants involved.

Let $f_n$ be a smooth approximation of $\mu$ in the $\ast$-weak topology of 
$\Mb$, such that 
$
||f_n||_{\Luno}\,\leq\,||\mu||_{\Mb}
$,
and let $u_n$ be the solutions of
$$
\cases{{\cal A}u_n=f_n\quad\hbox{ in }\Hduale\crr
	u_n\in\Hunozero.}
$$

Consider, for $\delta>0$, the Lipschitz continuous functions $h_\delta$ 
defined by
$$
\cases{h_\delta(s)=1\qquad\quad&if $|s|\leq k$\crr
	h_\delta(s)=0&if $|s|\geq k+\delta$\crr
	|h'_\delta(s)|={\displaystyle{1\over\delta}}&if 
	$k\leq|s|\leq k+\delta$,}
$$
and $S_\delta$ defined by
$$
\cases{S_\delta(s)=0\qquad\quad&if $|s|\leq k$\crr
	S_\delta(s)={\rm sign}(s)&if $|s|\geq k+\delta$\crr
	S'_\delta(s)={\displaystyle{1\over\delta}}&if $k\leq|s|\leq k+\delta$.}
$$

Using the equation, we can see that 
$
-{\rm div}(h_\delta(u_n) A(x)\nabla u_n)
$ belongs to $\Luno$ and that

$$
\eqalign{
&\intl_\Omega |-{\rm div}(h_\delta(u_n) A(x)\nabla u_n)|\,dx\crr
&\leq\intl_\Omega 
  |f_n|\left(h_\delta(u_n)+S_\delta^+(u_n)+S_\delta^-(u_n)\right)\,dx\crr
&=\intl_\Omega |f_n|\,dx\leq||\mu||_{\Mb}.}
$$
This implies 
$$
||\mu_k||_{\Mb}\,\leq\,||\mu||_{\Mb},
$$
(recall that $\mu_k={\cal A}T_k(\um)$) 
and we conclude 
thanks to Lemma \rf{piumeno}.
\endproof

\parag{\chp{finedimondo}}{The general existence theorem}

We come now to prove the existence and uniqueness of the solution to the 
Obstacle Problem, without the technical assumption that the 
obstacle be negative. From now on the only hypothesis will be (\rf{ipomin}).

\th{\thm{esist2}}{Let $\psi$ satisfy (\rf{ipomin}) and let $\mu\in\Mb$. 
Then there exists a (unique) solution of $OP(\mu,\psi)$.}

\proof{}
It is enough to show that we can refer to the case $\psi\leq 0$.
Indeed define
$$
\varphi:=\psi-u_\rho,
$$
which is, obviously, negative. 

By Theorem \rf{esistnegativo} there exists $v$ minimum in 
${\cal F}_\varphi(\mu-\rho)$, and
we prove that the function $u:=v+u_\rho$ is the minimum of $\Fpsimu$.

Trivially $u\geq\psi$ and, denoted the positive obstacle reaction
associated to $v$ by $\lambda$, we have
$u=v+u_\rho=\um+\ul$,
which shows that $u$ is an element of $\Fpsimu$.

Consider  now a function $w\in\Fpsimu$. By similar computations we deduce
that $w-u_\rho$ belongs to ${\cal F}_\varphi(\mu-\rho)$ and, by the
minimality of $v$,
$v\leq w-u_\rho$, so that we conclude $u\leq w\qe$,
and $\lambda$ is the obstacle reaction associated to $u$.
\endproof

\rem{\thm{Lewyvera}}{}From the previous proof we deduce that in the 
general case we have the inequality 
$$
||\lambda||_{\Mb}\,\leq\,||(\mu-\rho)^-||_{\Mb}.\eqno(\frm{mumenorho})
$$

\parag{\chp{nuovavita}}{Some stability results}

In this section we want to show some results of continuous dependence of
the solutions on the data.

The following proposition concerns the problem of stability with
respect to the obstacle, which, however, is not true in general (see Remark 
\rf{contro}).

\prop{\thm{psienne}}{Let $\psi_n:\Omega\to\overline{\R}$ be obstacles such
that
$$
\psi_n\leq\psi \quad\hbox{ and }\quad\psi_n\to\psi\quad \hbox{ q.e.   
in }\Omega,
$$
$\psi$ satisfies (\rf{ipomin}), and let
$u_n$ and $u$ be the solutions of $OP(\mu,\psi_n)$ and
$OP(\mu,\psi)$, respectively. Then
$$
u_n\to u\quad\hbox{ strongly in }\Wunoqzero.
$$
We also obtain that $u_n\to u\qe$ and that $T_k(u_n)\wto T_k(u)\weak2$,
for all $k>0$.}
\proof Since $u$ is trivially in ${\cal F}_{\psi_n}(\mu)$ for any $n$ we
have
$$
u_n \leq u\quad \qe.\eqno(\frm{inmezzo})
$$

To every minimum $u_n$ there corresponds a positive obstacle reaction
 $\lambda_n$,
satisfying inequality (\rf{mumenorho}), so we obtain that, up to a
subsequence,
$$
\eqalign{
\lambda_n\wto \hat\lambda&\quad\weaks\cr
u_n\to \hat u&\quad\hbox{ strongly in } \Wunoqzero\cr
} 
$$
and
$$
\hat u=\um+u_{\hat\lambda}.
$$
Hence, from (\rf{inmezzo}), $\hat u\leq u\ae$, and also q.e. On the
other side, we have to prove that
$\hat u\geq\psi$ q.e. in $\Omega$, in order to obtain $\hat u\in\Fpsimu$,
and so $u\leq \hat u$ q.e. in $\Omega$.

First consider the case when $\psi_n\leq\psi_{n+1}\qe$.

From this fact it follows that $u_n\leq u_{n+1}\qe$, and then $T_k(u_n)
\leq T_k(u_{n+1})$ q.e. in $\Omega$, for all $k>0$. Hence this sequence
has a quasi
everywhere limit.
On the other hand, the fact that $\mu+\lambda_n\wto\mu+\hat\lambda\weaks$
implies that
$ T_k(u_n)\wto T_k(\hat u)$ weakly in $\Hunozero$ and then, by Lemma 1.2
of [\rf{DAL-GAR}],
$T_k(u_n)\to T_k(\hat u)$ q.e. in $\Omega$. Since this holds for all $k>0$
we get also
$$
u_n\to \hat u\qe.
$$
Then, passing to the limit in $u_n\geq\psi_n\qe$ we get $\hat u\geq\psi
\qe$.
   
If the sequence $\psi_n$ is not increasing, consider
$$
\f_n:=\inf_{k\geq n} \psi_k,\eqno(\frm{cresce})
$$
so that $\f_n\nearrow \psi\qe$ and $\f_n\leq\psi_n\qe$. If $\overline u_n$
is the solution of $OP(\mu,\f_n)$ it is easy to see, using Definition
\rf{disvar}, that $\overline u_n\leq u_n \leq u\qe$. Applying the first
case to $\overline u_n$ and passing to the limit we get $u_n\to u\qe$.
\endproof

\bigskip

As for stability with respect to the right-hand side,
we will show later that in general it is not true that if 
$$
\mu_n\wto\mu\weaks
$$
then
$$
u_n\to u \strongq,
$$
where $u_n$ and $u$ are the solutions relative to $\mu_n$ and $\mu$ with
the fixed obstacle $\psi$.

However we can give now the following stability
result.

\prop{\thm{fortino}}{Let $\mu_n$ and $\mu$ be measures in $\Mb$ such that
$$
\mu_n\to\mu\hbox{ strongly in }\Mb,
$$
then 
$$
u_n\to u\hbox{ strongly in }\Wunoqzero
$$
where $u_n$ and $u$ are the solutions of $OP(\mu_n,\psi)$ and 
of $OP(\mu,\psi)$, respectively.}

\proof
Let $\lambda_n$ be the obstacle reactions associated to $u_n$, then
$$
\|\lambda_n\|_{\Mb}\,\leq\,\|(\mu_n-\rho)^-\|_{\Mb},
$$
so, up to a subsequence,
$$
\lambda_n\wto \hat\lambda\weaks
$$
and
$$
\eqalign{
u_n&\to \hat u\strongq\crr
T_k(u_n)&\wto T_h(\hat u)\weak2\ \forall k>0\cr}
$$
where $\hat u=\um+u_{\hat \lambda}$.
  
As  $T_k(u_n)\geq T_k(\psi)\qe$ for every $k\geq 0$, and for every $n$, we
have $T_k(\hat u)\geq T_k(\psi)$ q.e. in $\Omega$ for every $k>0$.

Passing
to the limit as $k\to +\infty$ we obtain that  $\hat u$ belongs
to $\Fpsimu$.

Let $v\in\Fpsimu$, with $\nu$ the associated measure.
Consider now $v_n$ the Stampacchia solution relative to
$\zeta_n:=\mu_n+(\mu_n-\mu)^-+\nu$.
Since $\zeta_n\to\mu+\nu$ strongly in $\Mb$, the sequence $v_n$
converges strongly in $\Wunoqzero$ to $v$.

Moreover
$v_n\geq v\geq \psi\qe$; hence $v_n\in{\cal F}_\psi(\mu_n)$, then $u_n\leq
v_n\qe$, and, in the limit,
$$
\hat u\leq v\ae,
$$ 
and hence also$\qe$.

\endproof

\rem{\thm{paragone}}{}
Thanks to this last result we can say that the solutions obtained in this
paper coincide with those given by Boccardo and Cirmi in [\rf{BOC-CIR}] 
and [\rf{BOC-CIR2}] when the data are $\Luno$-functions.

\bigskip

As said above we give now the counterexample showing that in general there
is not stability with respect to $\ast$-weakly convergent data.

\goodbreak\bigskip\noindent{\bf Example \thm{contro}} Let
$\Omega=(0,1)^N$ with $N\geq 3$, ${\cal A}=-\Delta$ and $\psi\equiv 0$.

The construction follows the one made by Cioranescu and Murat in
[\rf{CIO-MUR}].

For each $n\in {\rm I\! N}$, divide the whole of $\Omega$ into small
cubes  of side $\displaystyle{1\over n}$. In the centre of each of them
take  two balls: $B_{1\over{2n}}$, inscribed in the cube,  and
$B_{r_n}$ of ray $r_n=\displaystyle\left({1\over{2n}}\right)^{N\over{N-2}}$.

In each cube define $w_n$ to be the capacitary potential of $B_{r_n}$ with
respect to
$B_{1\over{2n}}$ extended by zero in the rest of the cube.

Hence 
$$
\Delta w_n=\mu_n,
$$
with
$$
\mu_n\wto 0 \hbox{ \ both weakly in }\Hduale\hbox{ and}\weaks.
$$
(see [\rf{CIO-MUR}]). Thus $ w_n\wto 0$ weakly in $\Hunozero$.

Let $u_n$ be the solution of $VI(\mu_n,0)$. Using $w_n$ as test function 
in the Variational Inequality we get $||u_n||_{\Hunozero}\leq C$. By 
contradiction assume that its $\Hunozero$-weak limit is zero.

Consider the function $z_n:=u_n+w_n$ which must then converge to zero
weakly in $\Hunozero$. Obviously $z_n\geq w_n$ q.e. in $\Omega$ and then
$z_n\geq 1$ on $\bigcup B_{r_n}$.
Hence if we define the obstacles
$$
\psi_n:=\left\{\eqalign{&1\quad\hbox{ in }\bigcup B_{r_n}\crr
			&0\quad\hbox{ elsewhere}\cr}\right. 
$$
$z_n\geq\psi_n$. Call $v_n$ the function realizing 
$$
\min_{{v\geq\psi_n}\atop{v\in\Hunozero}}\intl_\Omega|\nabla v|^2\,dx.
$$

A simple computation yields
$$
-\Delta z_n\,=\,-\Delta u_n - \Delta w_n\,\geq\, 0.
$$

Then $z_n\geq v_n\geq 0$, so that
$$
v_n\wto 0\hbox{ weakly in } \Hunozero.
$$
But this is not possible because a $\Gamma$-convergence result contained
in [\rf{CAR-COL}] says that there exists a constant $c>0$ such that
$v_n$ tends to the minimum point of
$$
\min_{{v\geq 0}\atop{v\in\Hunozero}}\intl_\Omega|\nabla v|^2\,dx+
c\intl_\Omega|(v-1)^-|^2\,dx
$$
which is not zero.

\parag{\chp{ipotesi}}{Comparison with the classical solutions}

As announced, in this section, we want to show that the new formulation 
of Obstacle Problem is consistent with the classical one.

To talk about the equivalence of the two formulations it is necessary 
that both make sense. So we will work under the hypothesis that
$\mu\in\Inters$ and that the 
obstacle $\psi$ satisfies

$$
\exists z\in\Hunozero \hbox{ s.t. }z\geq\psi\qe;\eqno(\frm{ipo1})
$$
$$
\exists\rho\in\Mbp \hbox{ s.t. }u_\rho\geq\psi\qe.\eqno(\frm{ipo2})
$$

Later on we will discuss these conditions in deeper details.

\lemma{\thm{equiv2}}{If there exists a measure $\sigma\in\Inters$ such that 
$u_\sigma\geq\psi\qe$, then the 
solutions of $VI(\mu,\psi)$ and of $OP(\mu,\psi)$ coincide.}
\proof
Let $u$ be the solution of $VI(\mu,\psi)$. 
Subtracting $u_\sigma$ to it, and with the same technique as in the proof of 
Theorem \rf{esist2}, we return to the case of negative obstacle and we can 
use Lemma \rf{equiv1}. \endproof

\th{\thm{equiv3}}{Under the hypotheses (\rf{ipo1}) and (\rf{ipo2}), 
the solutions of $VI(\mu,\psi)$ and of $OP(\mu,\psi)$ coincide.}

\proof As a first step consider the case of an obstacle bounded from 
above by a constant $M$.
The measure $\rho_M:={\cal A}(T_M(u_\rho))$ is in $\Inters$ and 
$T_M(u_\rho)\geq\psi$ so that we are in the hypotheses of the previous lemma.

If, instead, $\psi$ is not bounded, we consider $\psi\wedge k$, and, with
respect to this new obstacle, conditions (\rf{ipo1}) and 
(\rf{ipo2}) are satisfied by the function $T_k(u_\rho)$. 

Hence we can apply the first step and say that $u_k$, solution of 
$VI(\mu,T_k(\psi))$, is also the solution of $OP(\mu,T_k(\psi))$.

From the classical theory we know that the sequence $u_k$ tends in 
$\Hunozero$ to the solution of $VI(\mu,\psi)$, while from Proposition 
\rf{psienne} 
$u_k$ converges in $\Wunoqzero$ to the solution of $OP(\mu,\psi)$.
\endproof

\medskip

A little attention is required in treating conditions (\rf{ipo1}) and 
(\rf{ipo2}). Each 
one is necessary for the corresponding problem to be nonempty, but together 
they can be somewhat weakened.

First of all we underline that no one of the two conditions is implied by 
the other. This is seen with the following examples.

\ex{\thm{1no2}}{}Let $\Omega=(-1,1)\subset\R$ and let ${\cal 
A}=-\Delta=-u''$. Take $\psi\in {\rm H}^1_0(-1,1)$ such that $-\psi''$
is an unbounded positive Radon measure. For instance we may take
$\psi=(1-|x|)(1-\log (1-|x|))$.

Now (\rf{ipo1}) is trivially true, 
and the solution of $VI(0,\psi)$ is $\psi$ itself.
If also (\rf{ipo2}) were true, then $\psi$ would be also the solution
of $OP(0,\psi)$. But this is not possible, because,
being $-\psi''$ an unbounded measure, we can not write it as $\ul$ 
for some $\lambda\in\Mbp$.

\ex{\thm{2no1}}{}Let $N\geq 3$, ${\cal A}=-\Delta$ and $\rho=\delta_{x_0}$, 
the Dirac delta 
in a fixed point $x_0\in\Omega$.

Take $\psi=u_{\delta_{x_0}}$, the Green function with pole at $x_0$. 
Then (\rf{ipo2}) holds, but if also (\rf{ipo1}) 
held we would have $\psi\in{\rm L}^{2^*}(\Omega)$ which is not true.

\medskip 

On the other side we already saw in the proof of Theorem \rf{equiv3} that 
if we add to Condition (\rf{ipo2}) the assumption that the obstacle be bounded, 
this is enough for (\rf{ipo1}) too to hold.

Moreover, if, besides (\rf{ipo1}), we assume that the obstacle is ``controlled near 
the boundary'' also Condition (\rf{ipo2}) is true:

Assume that (\rf{ipo1}) holds and there exists a compact 
$J\subset\Omega$, such that $\psi\leq 0$ in $\Omega\setminus J$. Then also 
(\rf{ipo2}) holds.
Indeed just take as $\rho$ the obstacle reaction corresponding
 to $u$, the solution of 
$VI(0,\psi)$. Then 
$$
{\rm supp}\rho\subset J,
$$
and $\rho\in\Mbp$.

A finer condition expressing the ``control near the boundary'' is 
\item{(2')}$
\exists J\hbox{ compact }\subset\Omega\hbox{ and }\exists 
\tau\in\Mbp\cap\Hduale\ :\ u_{\tau}\geq\psi\hbox{ in }\Omega\setminus J.$

\bigskip

In conclusion we want to remark that, in general, in classical Variational 
Inequalities, the obstacle reaction associated to the 
solution is indeed a Radon measure, but it is not always 
bounded, as Example \rf{1no2} shows.

On the other side, in the new setting, the minimum of $\Fpsimu$ is not, 
in general, an element of $\Hunozero$.

Hence the two formulations do not overlap completely and no one is 
included in the other.

\parag{\chp{approx}}{Approximation properties}

As we have seen so far, if we have a 
sequence $\mu_n$ $\ast$-weakly convergent to $\mu$, we can not deduce 
convergence of solutions, but, from (\rf{mumenorho}) we have
$$
||\lambda_n||_{\Mb}\leq||(\mu-\rho)^-||_{\Mb},
$$
where the $\lambda_n$ are the obstacle reactions relative to 
the solutions $u_n$.
So, up to a subsequence, 
$$
\lambda_n\wto\hat\lambda\weaks
$$
and
$$
u_n\wto\hat u\,=\um+u_{\hat\lambda}\strongq.
$$
With the same argument used in the proof of
Theorem~\rf{esistnegativo} we can show that $\hat u\geq \psi\qe$. Hence
$\hat u\geq u$, the minimum of $\Fpsimu$. 

On the other hand, in Theorem \rf{esist2} we have obtained the 
solution of $OP(\mu,\psi)$
as a limit of the solutions to $OP({\cal A}T_n(u_{\mu-\rho})+\rho,\psi)$. 
We remark that if $\rho$ belongs to the ordered dual of $\Hunozero$ that 
is $V:=\left\{\mu\in
\Inters\ :\ |\mu|\in\Hduale\right\}$, then 
the approximating problems are actually Variational Inequalities.

Thanks to these two facts we can 
characterize the solution $u$ of $OP(\mu,\psi)$ by 
approximation with solutions of Variational Inequalities with data in
$V$ as follows.

\medskip

\item{1.} For every sequence $\mu_n$ in $\Mb$, with
 $\mu_n\wto\mu\weaks$, we have
$$
s\hbox{-}\Wunoqzero\hbox{-}\lim_{n\to\infty} u_n\geq u.
$$

\item{2.} There exists a sequence $\mu_n\in V$, with
 $\mu_n\wto\mu\weaks$ such that 
$$
s\hbox{-}\Wunoqzero\hbox{-}\lim_{n\to\infty} u_n= u
$$

In other words:
$$
u\,=\,\min\left\{ s\hbox{-}{\hskip-.3cm}\lim_{n\to+\infty}u_n\ :\ u_n 
\hbox{ sol. }VI(\mu_n,\psi) ,\,\mu_n\in V,\,\mu_n\wto\mu\weaks\right\}. 
$$

\parag{\chp{emmezero}}{Measures vanishing on sets of zero Capacity}

We show now an example (suggested by L. Orsina and A. Prignet) in 
which the solution of the Obstacle Problem with right-hand side 
measure does not touch the obstacle, though it is not the solution of the 
equation.

\ex{\thm{delta}}{Let $N\geq 2$, $\Omega$ be the ball $B_1(0)$, and ${\cal 
A}=-\Delta$. Take the datum $\mu$ a negative measure concentrated on a 
set of zero 2-Capacity and the obstacle $\psi$ negative and bounded below
by a
constant $-h$.
Let $u$ be the solution of $OP(\mu,\psi)$, then 
$u=\um+\ul$. We want to show that $\lambda=-\mu$.

First observe that, for minimality, $u\leq0$; on the other hand $u\geq 
-h$, so that $u=T_h(u)$ and hence $u\in\Hunozero$. This implies that the 
measure $\mu+\lambda$ is in $\Inters$, which is contained in $\Mbo$, the 
measures which are zero on the sets of zero 2-Capacity (see 
[\rf{BOC-GAL-ORS}]). In other words $\lambda=-\mu+{\hat \lambda}$, with 
${\hat \lambda}$ a measure in $\Mbo$, and so positive, since $\lambda$ is 
positive. Then $u\geq 0$, and finally $u=0$.   
Thus the solution can be far above the obstacle, but the obstacle reaction
is 
nonzero, and is exactly $-\mu$.}

\bigskip

\rem{\thm{contro}}{This example shows also that in general there is no 
continuous dependence on the obstacles. Indeed, if $h\to +\infty$, then 
the solution of $OP(\mu,-h)$ is identically zero for each $h$, while the 
solution of $OP(\mu,-\infty)$ is $\um$.}
\bigskip
 
We want to 
consider here a class of data for which the above phenomenon is avoided.

Consider, as datum, a measure in $\Mbo$. In this case we can use 
the fact (contained in [\rf{BOC-GAL-ORS}]) that for any such 
measure $\mu$ 
there exists a  function $f$ in $\Luno$ and a functional $F$ in $\Inters$, 
such that $\mu=f+F$. If, in addition $\mu\geq0$, then also $f$ can be 
taken to be positive.

We want to show that also the obstacle reaction $\lambda$ belongs to 
$\Mbo$ and that in this particular case we can write our Obstacle Problem
in a variational way, that is with and ``entropy formulation".

We begin by considering the case of a negative obstacle.

\lemma{\thm{reazioni}}{Let $\psi\leq0$ and let $\mu_1,\,\mu_2\in\Inters$. 
Let $\lambda_1$ and $\lambda_2$ be the reactions of the obstacle 
corresponding to the solutions $u_1$ and $u_2$ of $VI(\mu_1,\psi)$ and 
$VI(\mu_2,\psi)$, 
respectively. 

If $\mu_1\leq\mu_2$ then $\lambda_1\geq\lambda_2$.}
\proof
This proof is inspired by Lemma 2.5 in [\rf{DAL-SKR}]. We easily have that 
$u_1\leq u_2$.

Take now a function $\varphi\in {\cal D}(\Omega)$, $\varphi\geq0$, 
and set $$
\fe:=\e\f\wedge(u_2-u_1)\,\in\Hunozero.
$$
Now, using the hypothesis that $\mu_1\leq\mu_2$ and monotonicity of ${\cal 
A}$, compute $$
\eqalign{
\dualita{\lambda_1}{\e\f-\fe}&\geq\dualita{{\cal A}u_1}{\e\f-\fe}-
                                \dualita{\mu_2}{\e\f-\fe}\crr
	&=\dualita{{\cal A}u_1-{\cal A}u_2}{\e\f-\fe}
		+\dualita{\lambda_2}{\e\f-\fe}\crr 
	&\geq\e{\hskip-.7cm}\intl_{\left\{u_2-u_1\leq\e\f\right\}}
	{\hskip-.7cm} A(x)\nabla(u_1-u_2)\nabla\f
	+\e\dualita{\lambda_2}{\f}-\dualita{\lambda_2}{\fe}.\cr}
$$
Now, using $u_1$ as a test function in $VI(\mu_2,\psi)$ and the fact that 
$u_2-u_1\geq\fe\geq0$ we easily get $\dualita{\lambda_2}{\fe}=0$.

Since, also, $-\dualita{\lambda_1}{\fe}\leq0$ we obtain
$$
\dualita{\lambda_1}{\f}\geq{\hskip-.7cm}\intl_{\left\{u_2-u_1\leq\e\f\right\}}
{\hskip-.7cm} A(x)\nabla(u_1-u_2)\nabla\f+\dualita{\lambda_2}{\f}.
$$
Passing to the limit as $\e\to 0$ and observing that 
$$
{\hskip-.7cm}\intl_{\left\{u_2-u_1\leq\e\f\right\}}     
        {\hskip-.7cm} A(x)\nabla(u_1-u_2)\nabla\f\longrightarrow 
{\hskip-.4cm}\intl_{\left\{u_2=u_1\right\}}     
        {\hskip-.4cm} A(x)\nabla(u_1-u_2)\nabla\f=0,
$$
we get the thesis.
\endproof

Let us see now what can we say more if $\mu\in\Mbo$, still in the case of 
negative obstacle.

\lemma{\thm{lambdazero}}{Let $\psi\leq 0$ and let $\mu\in\Mbo$ then the 
obstacle reaction relative to the solution of $OP(\mu,\psi)$ is also in 
$\Mbo$.}

\proof
It is not restrictive to assume $\mu$ to be negative. Indeed, if 
$\mu=\mu^+-\mu^-$, then also $\mu^+$ and $\mu^-$ are in $\Mbo$. Hence the 
minimum of $\Fpsimu$ can be written as $u_{\mu^+}+v$ with $v$ minimum in 
${\cal F}_{\psi-u_{\mu^+}}(-\mu^-)$, and the same obstacle reaction 
$\lambda$; and so we are in the case of a negative measure.

Consider now the decomposition $\mu=f+F$ with $f\leq 0$. And let 
$\mu_k:=T_k(f)+F$ so that $\mu_k\to\mu$ strongly in $\Mb$.

Let $u_k$ be the solution of $OP(\mu_k,\psi)$. It is also the solution of 
$VI(\mu_k,\psi)$ so that $\lambda_k\in\Mbo$.

Thanks to Proposition \rf{fortino} we have that $u_k\to u=\um+\ul$ 
strongly in $\Wunoqzero$ and that $\lambda_k\wto\lambda\weaks$.

From the fact that $\mu_k\geq\mu_{k+1}$ and from Lemma \rf{reazioni} we 
obtain that $\lambda_k\leq\lambda_{k+1}$. Hence if we define 
$$
\hat\lambda(B):=\lim_{k\to\infty}\lambda_k(B)\qquad\forall B\hbox{ Borel 
set in }\Omega,
$$
we know from classical measure theory that it is a bounded Radon measure, it 
is in 
$\Mbo$, since all $\lambda_k$ are, and necessarily coincides with $\lambda$. 
So $\lambda\in\Mbo$.
\endproof
 
\bigskip

In order to pass to a signed obstacle observe first that the minimal 
hypothesis (\rf{ipomin}) becomes necessarily 
$$
\exists\sigma\in\Mbo\,:\,u_\sigma\geq\psi.\eqno(\frm{ipo3})
$$
Once we have noticed this, it is easy to use the result for a negative 
obstacle, as we did in the proof of Theorem \rf{esist2} and obtain the 
following result.

\th{\thm{lambda}}{Let $\psi$ satisfy hypothesis (\rf{ipo3}), and let 
$\mu$ be in $\Mbo$. Then the obstacle reaction relative to the solution 
of $OP(\mu,\psi)$ belongs to $\Mbo$ as well.}

\bigskip

\rem{\thm{forte}}{Notice that thanks to the pointwise convergence we have, 
in this case, that $\lambda_k\to\lambda$ strongly in $\Mb$.}

\bigskip

\rem{\thm{entropia}}{These properties of the case of $\Mbo$ measures, 
allow us to write the Obstacle Problem in a ``more variational'' way. 
Namely, if $\mu\in\Mbo$ and its decomposition is $\mu=f+F$ then the 
function $u$ solution of $OP(\mu,\psi)$ satisfies also
$$
\cases{
	\dualita{{\cal A}u}{T_j(v-u)}\geq{\displaystyle\intl_\Omega} 
		f\,T_j(v-u)+\dualita{F}{T_j(v-u)}\crr
	\forall v\in\Hunozero\cap\Linf,\ v\geq\psi\qe\cr}.
$$
This is similar to the entropy formulation given by Boccardo and Cirmi in 
[\rf{BOC-CIR}] in the case of datum in $\Luno$. The proof that such a 
formulation holds is made by 
approximation.
To this aim we choose a particular sequence of measures 
$\mu_k:=T_k(f)+F$, so that $\mu_k\to\mu$ strongly in $\Mb$. Hence also the 
solutions of $OP(\mu_k,\psi)$ (and also of $VI(\mu_k,\psi)$) converge 
strongly in $\Wunoqzero$ to $u$ solution of $OP(\mu,\psi)$. Then $u_k$ solves
$$
\cases{
	\dualita{{\cal A}u_k}{v-u_k}\geq\dualita{\mu_k}{v-u}\crr
	\forall v\in\Hunozero,\ v\geq\psi\qe\cr}
$$
In this inequality we can use as test functions $v=T_j(w-u_k)+u_k$, with 
$w\in \Hunozero \cap \Linf$, $w\geq\psi\qe$, 
and, by calculations similar to those in [\rf{BOC-CIR}], get the result.}

\intro{References}
\def\interrefspace{\smallskip}  
{\ninepoint\frenchspacing


\item{[\bib{BOC-CIR}]}BOCCARDO  L., CIRMI  G.R.: Nonsmooth unilateral
problems. In {\it Nonsmooth optimization: methods and applications}
(Proceedings Erice 1991) Eds Giannessi, Gordon and Breach, (1992) 1-10.
\interrefspace

\item{[\bib{BOC-CIR2}]}BOCCARDO  L., CIRMI  G.R.: Existence and uniqueness
of solutions of unilateral problems with $\rm L^1$-data. To appear.
\interrefspace

\item{[\bib{BOC-GAL}]}BOCCARDO  L., GALLOU\"ET  T.: Probl\`emes
unilat\'eraux avec donn\'ees dans $\rm L^1$. 
{\it C. R. Acad. Sci. Paris, S\'erie I \/} {\bf 311} (1990), 617-619. 
\interrefspace

\item{[\bib{BOC-GAL-ORS}]}BOCCARDO  L., GALLOU\"ET  T., ORSINA  L.: 
Existence and uniqueness of entropy solutions for nonlinear elliptic 
equations with measure data. {\it Ann. Inst. Henri Poincar\'e}  {\bf 13} 
(1996), 539-551.
\interrefspace 

\item{[\bib{BOC-MUR}]}BOCCARDO  L., MURAT  F.: A property of nonlinear
elliptic equations when the right-hand side is a measure. 
{\it Potential Analysis \/} {\bf 3} (1994), 257-263. 
\interrefspace

\item{[\bib{CAR-COL}]}CARBONE  L., COLOMBINI  F.: On convergence of
functionals with unilateral constraints. 
{\it J. Math. Pures et Appl. \/} {\bf 59} (1980), 465-500. 
\interrefspace

\item{[\bib{CIO-MUR}]}CIORANESCU  D., MURAT  F.: Un terme \'etrange venu
d'ailleurs I and II. In  
{\it Nonlinear partial differential equations and their applications. \/}
Pitman research notes in mathematics; n. 60, 98-138; n. 70, 154-178.
\interrefspace


\item{[\bib{DAL-4}]}
DAL MASO  G.:
Some necessary and sufficient conditions for the convergence of sequences
of unilateral convex sets
{\it J. Funct. Anal.} {\bf 62}, 2 (1985), 119-159.
\interrefspace

\item{[\bib{DAL-GAR}]}  
DAL MASO  G., GARRONI A.:   
The capacity method for asymptotic Dirichlet problems.
{\it Asymptotic Anal.} {\bf 15} (1997), 299-324.
\interrefspace

\item{[\bib{DM-MAL}]}DAL MASO  G., MALUSA  A.: Some properties of reachable 
solutions of nonlinear elliptic equations with measure data.
Preprint S.I.S.S.A., April 1997.
\interrefspace

\item{[\bib{DAL-SKR}]}DAL MASO  G., SKRYPNIK  I.V.:
Capacity theory for monotone operators.
Preprint S.I.S.S.A., January 1995.
\interrefspace


\item{[\bib{HEI-KIL}]}HEINONEN  J., KILPEL\"AINEN  T., MARTIO  O.: {\it
Nonlinear potential theory of degenerate elliptic equations.} Clarendon
Press, Oxford, 1993.
\interrefspace

\item{[\bib{KIL-MAL}]}KILPEL\"AINEN  T., MAL\'Y  J.: Degenerate elliptic
equations with measure data and nonlinear potentials. 
{\it Ann. Scuola Norm. Sup. Pisa \/} {\bf 19} (1992), 591-613. 
\interrefspace

\item{[\bib{KIN-STA}]}KINDERLEHRER  D., STAMPACCHIA  G.: {\it
An introduction to Variational Inequalities and their applications.} 
Academic, New York, 1980.

\item{[\bib{PRI}]}PRIGNET  A.:
Remarks on existence and uniqueness of solutions of elliptic problems 
with right-hand side measures. {\it Rend. Mat.} {\bf 15} (1995), 97-107.
\interrefspace

\item{[\bib{SER}]}SERRIN  J.: Pathological solutions of elliptic
differential equations. {\it Ann. Scuola Norm. Sup. Pisa} {\bf 18} 
(1964), 385-387.
\interrefspace

\item{[\bib{STA}]}STAMPACCHIA  G.: Le probl\`eme de Dirichlet pour les
\'equations elliptiques du second ordre \`a coefficients discontinus. 
{\it Ann. Inst. Fourier Grenoble \/} {\bf 15} (1965), 189-258. 
\interrefspace

\item{[\bib{TRO}]}TROIANIELLO  G.M.: {\it
Elliptic differential equations and obstacle problems.} 
Plenum Press, New York, 1987.

\bye